\documentclass[11pt,dvipdfmx]{article}

\usepackage[top=30truemm,bottom=30truemm,left=25truemm,right=25truemm]{geometry}

\usepackage{amsfonts}
\usepackage{amssymb}
\usepackage{amsmath}
\usepackage{amsthm}

\DeclareMathOperator{\pfaff}{pf}

\DeclareMathOperator{\mat}{Mat}
\DeclareMathOperator{\sgn}{sgn}
\DeclareMathOperator{\head}{head}
\DeclareMathOperator{\fd}{fd}
\DeclareMathOperator{\bd}{bd}
\DeclareMathOperator{\Padj}{Padj}

\newtheorem{theorem}{Theorem}
\newtheorem{lemma}{Lemma}
\newtheorem{definition}{Definition}
\newtheorem{corollary}{Corollary}

\newcommand{\vnc}[1]{\mathbf{VNC}^{#1}}

\newcommand{\vsl}{\mathbf{V\#L}}

\newcommand{\sharpl}{\#\mathsf{L}}
\newcommand{\ext}{\exists}
\newcommand{\all}{\forall}
\newcommand{\sbi}[1]{\Sigma^B_{#1}}
\newcommand{\pairs}[1]{\langle #1\rangle}

\title{Formalizing Pfaffian in bounded arithmetic}
\author{Satoru Kuroda\footnote{This work was supported by JSPS KAKENHI Grant Number 18K03400.}
\\Department of Culture and Informatics,\\
Gumma Prefectural Women's University}

\date{}

\begin{document}

\maketitle

\section{Introduction}

Proving theorems of linear algebra, especially properties of the determinant 
is a central theme in bounded reverse mathematics. Although the determinant 
has several $\sharpl$ algorithms, many of its properties are known to be provable 
in a theory slightly stronger than $\sharpl$. 

The seminal work of Soltys and Cook on bounded reverse mathematics of linear algebra 
revealed that some important properties such as the cofactor expansion, 
the axiomatic definition of the determinant and Cayley-Hamilton Theorem 
are equivalent over the theory $\vsl$. Also they proved that the multiplicativity 
of the determinant implies all these properties. 

Soon after, the celebrated result by Tzameret and Cook gave an upper bound on 
the provability of these properties. It is proved in \cite{t-c} that 
the multiplicativity of the determinant is provable in $\vnc{2}$. 

Also, Ken and the author \cite{k-k} showed that properties of matrix rank are provable 
in $\vnc{2}$ by using the result in \cite{t-c} and establishing the interpretation of 
extensions of Soltys theory for linear algebra in $\vnc{2}$. 
However, it is still open that the above properties of the determinant and matrix rank are 
provable in some weaker theories such as $\vsl$. In particular, the proof in \cite{t-c} 
is based on the algorithm for the determinant via Schur complement. On the other hand, 
faster algorithms such as Berkowitz algorithm \cite{berk} are formalizable in $\vsl$. 

In this article, we propose to extend the study of proof complexity of linear algebra 
along this line to Pfaffian. 

Pfaffian was introduced by Pfaff in 19th century in relation with partial differential equations. 
Recently, many applications are given in combinatorics and representation theory. 
Computing Pfaffian is very similar to computing the determinant and 
many fast algorithms for the determinant are generalized to Pfaffian 
which include the characterization via clow sequences. 

In this article, we will give a Berkowitz type algorithm for Pfaffian and 
prove its correctness by way of clow sequences technique which was developed by 
Mahajan, Subramanya and Vinay. This is used to formalize Pfaffian in the theory 
$\vsl$.

Then we also consider the provability of Pfaffian properties over the theory $\vsl$. 
Especially, we consider the problem of proving properties from Pfaffian version of 
multiplicativity. 

We also present a version of Cayley-Hamilton type theorem for Pfaffian. 
Cayley-Hamilton type theorem for Pfaffian has been unfamiliar until recently. 
By examining the proof of Cayley-Hamilton Theorem from cofactor expansion in 
\cite{c-s} carefully, we present a theorem which is equivalent to cofactor expansion 
and the axiomatic definition of Pfaffian. 

To author's knowledge, our version of Pfaffian Cayley-Hamilton is new and 
we expect that it can be used to prove various properties of linear algebra. 

Our goal is not only to extend the proof complexity problem of linear algebra 
but also to extend bounded reverse mathematics to combinatorics and representation theory. 
The final section is devoted to showing our perspective.

\section{Preliminaries}

Due to the space limit, we refrain from giving details of bounded arithmetic and complexity theory. 
We work in two sort bounded arithmetic developed by Cook and Ngyuen \cite{c-n}. 
The theory $\vsl$ consists of axioms 
\begin{itemize}
\item $\sbi{0}$-COMP:
\begin{equation*}
\all a\;\all x<a\;\ext Y\;(x\in Y\leftrightarrow Y(x))
\end{equation*}
where $\varphi(x)\in\sbi{0}$ does not contain $Y$. 
\item String Multiplication:
\begin{equation*}
\all X,Y\;\ext Z\;(Z=X\cdot Y)
\end{equation*}
\item Matrix Powering:
\begin{equation*}
\all X\mbox{ : square matrix}\;\all n\;\ext Y\;(Y=X^n)
\end{equation*}
\end{itemize}

The complexity class $\sharpl$ consists of functions which are logspace reducible to 
the determinant. It is known that matrix powering is complete for $\sharpl$ 
and thus we have

\begin{theorem}
A function is $\sbi{1}$ definable in $\vsl$ if and only if it is in $\sharpl$. 
\end{theorem}

Note that induction for $\sbi{0}$ formula is provable in $\vsl$ even when we 
extend the language by $\sbi{1}$ definable functions. This fact will be a crucial 
tool in proving matrix properties. 

Pfaffian is defined in a similar manner as for the determinant. 
Specifically, let $A\in\mat(2n,2n)$ be skew symmetric. Then its Pfaffian is defined as 
\begin{equation}
\pfaff(A)=\sum_{\sigma\in\mathcal{M}_{2n}}\sgn(\sigma)a_{\sigma(1)\sigma(2)}\cdots a_{\sigma(2n-1)\sigma(2n)}
\end{equation}
where $\mathcal{M}_{2n}$ represents the set of perfect matchings on $[2n]$ such that 
$$
\sigma(1)<\sigma(3)<\cdots\sigma(2n-1).
$$

Pfaffian can be regarded as a generalization of the determinant in the sense that 
$\det(A)$ for $n\times n$ matrix $A$ is computed by Pfaffian as 
\begin{equation}\label{eq:det-pf}
\det(A)=
(-1)^{n(n-1)}\pfaff
\begin{pmatrix}
0&A\\
-^t\!A&0
\end{pmatrix}
\end{equation}

For skew symmetric matrix $A\in\mat(2n,2n)$, the following relation is known:
\begin{theorem}[Cayley]
If $A\in\mat(2n,2n)$ is a skew symmetric matrix then 
\begin{equation}\label{eq:cayley}
\det(A)=\pfaff(A)^2
\end{equation}
\end{theorem}

Our formalization of Pfaffian is based on the characterization by way of clow sequences 
due to Mahajan, Vinay. A clow (closed walk) on $[n]$ is a list of edges 
$$
(i_1,i_2),(i_2,i_3),\ldots,(i_m,i_1)
$$
such that $i_1<i_k$ for all $2\leq k\leq m$. The first index $i_1$ is called the head of $C$ and 
is denoted by $\head(C)$. 

A pclaw is a list $E_1,E_2,\ldots,E_m$ 
where each $E_k$ is a pair $(e^k_1,e^k_2)$ of edges such that either
\begin{itemize}
\item $e^k_1=(i,2j-1)$ and $e^k_2=(2j-1,2j)$ or
\item $e^k_1=(i,2j)$ and $e^k_2=(2j,2j-1)$.
\end{itemize}

Let $C$ be a pclow. Define
$$
\fd(C)=\#\{(i,j)\in C : i<j\},\ \bd(C)=\#\{(i,j)\in C : i>j\}.
$$
and
$$
\sgn(C)=(-1)^{f(C)+1}. 
$$
For $A=(a_{ij})\in\mat(2n,2n)$, we define 
$a_{ij}^+=a_{ij}$ if $i<j$ and $a_{ij}^+=a_{ji}$ if $i>j$. 
The weight of a clow $C=\pairs{e_1,e_2,\ldots,e_{2m}}$ over $A$ is the product
$$
w_A(C)=\prod_{1\leq k\leq m}a_{e_{2k-1}}^+.
$$

A pclow sequence is a sequence $\bar{C}=\pairs{C_1,\ldots,C_l}$ of pclows such that 
$$
\head(C_1)=1<\head(C_2)<\cdots<\head(C_l).
$$
We define the sign and the weight of a pclow sequence as 
$$
\sgn(\bar{C})=\prod_{C\in\bar{C}}\sgn(C)\mbox{ and }
w_A(\bar{C})=\prod_{C\in\bar{C}}w_A(C)
$$
respectively. Finally the length of a pclow or a pclow sequence is the number of edges 
occuring in it. 

\begin{theorem}[Mahajan et.al.]
Let $A\in\mat(2n,2n)$ be skew symmetric. Then 
$$
\pfaff(A)=\sum_{\bar{C}\mbox{:pclow seq. }|\bar{C}|=2n}\sgn(\bar{C}))w_A(\bar{C}).
$$
\end{theorem}

\section{Berkowitz-type algorithm for Pfaffian}

In this section we construct a $\sharpl$ algorithm for Pfaffian. 

For $n\in\omega$, we define the skew symmetric matrix $J_n\in\mat(2n,2n)$ by 
$$
J_1=
\begin{pmatrix}
0&1\\
-1&0
\end{pmatrix},\ 
J_n=
\left.
\begin{pmatrix}
J_1&&&\\
&J_1&&\\
&&\ddots&\\
&&&J_1
\end{pmatrix}
\right\}n\mbox{ times}\ (n\geq 1).
$$
We omit the subscript if it is clear from the context. 

\begin{definition}[PB algorithm]\label{def:pb algorithm}
Let $A\in\mat(2n,2n)$ be skew symmetric and 
$$
\begin{pmatrix}
0&a_{12}&R\\
-a_{12}&0&-^tS\\
-^tR&S&M
\end{pmatrix}
$$
be its block decomposition. Define Berkowitz algorithm $P_A$ as 
\begin{equation}
P_A=
\begin{pmatrix}
1&&&\\
a_{12}&1&&&\\
RJS&a_{12}&\ddots&&\\
RJ(MJ)S&RJS&\ddots&\ddots&\\
\vdots&\ddots&\ddots&\ddots&\\
\vdots&\ddots&\ddots&\ddots&1\\
RJ(MJ)^{n-2}S&RJ(MJ)^{n-3}S&\cdots&\cdots&a_{12}
\end{pmatrix}
\in\mat(n+1,n).
\end{equation}
We define Pfaffian coefficients $\bar{P}_A=(p_n,p_{n-1},\ldots,p_0)$ as $\bar{P}_A=(1,a_{12})$ 
if $n=2$ and
$$
\bar{p}_A=P_A\bar{p}_M. 
$$
if $n>2$. 
\end{definition}

This algorithm is already suggested by Rote \cite{rote} in somewhat awkward manner. 
We present it here in a complete form and prove its correctness below.

We will show that PB algorithm computes Pfaffian. More generally we have 
\begin{theorem}\label{thm:pb}
Let $A\in\mat(2n,2n)$ be skew symmetric and $\bar{P}_A=(p_n,p_{n-1},\ldots,p_0)$ be its 
Pfaffian sequence. Then 
\begin{equation}
p_{n-k}=\sum_{\substack{\bar{C}\mbox{ : pclow seq. }\\|\bar{C}|=2k}}\sgn(\bar{C})w_A(\bar{C})
+\sum_{\substack{\bar{C}\mbox{ : pclow seq. on }[3,2n]\\|\bar{C}|=2k}}\sgn(\bar{C})w_M(\bar{C}).
\end{equation}
for $1\leq k\leq n-1$ and 
\begin{equation}
p_0=\sum_{\substack{\bar{C}\mbox{ : pclow seq. }\\|\bar{C}|=2n}}\sgn(\bar{C})w_A(\bar{C})
\end{equation}
Hence $\pfaff(A)=p_0$. 
\end{theorem}

To prove Theorem \ref{thm:pb}, we first notice that each entry in the matrix $P_A$ computes 
the sum of signed weights of clows. For instance, consider the entry $a_{12}$. 
The only possible clow starting from $(1,2)$ is $C=\pairs{(1,2),(2,1)}$ with $w_A(C)=a_{12}$. 
Moreover, note that $f(C)=1$ and thus $\sgn(C)=(-1)^{1+1}=1$. Hence we have
$$
\sum_{C\mbox{ : pclow }|\bar{C}|=2}\sgn(C)w_A(C)=a_{12}.
$$

In general we have 
\begin{lemma}\label{lemma:plow}
Let $A\in\mat(2n,2n)$ be skew symmetric with its block decomposition given as above. 
Then 
\begin{equation}\label{eq:pclow}
RJ(MJ)^{k-2}S=\sum_{\substack{C\mbox{ : pclow }\\|\bar{C}|=2k,\head(C)=1}}\sgn(C)w_A(C). 
\end{equation}
\end{lemma}

\begin{proof}
Let $A=(a_{ij})$ and $J=(b_{ij})$. 
Note that $RJ(MJ)^{k-2}S$ is the sum of the products of the form
$$
a_{e_1}b_{e_2}\cdots a_{e_{2k-1}}\cdot(-1).
$$
where $C=\pairs{e_1,e_2,\ldots,e_{2k-1},(2,1)}$ is a pclow. So we have 
$$
\begin{array}{rcl}
a_{e_1}b_{e_2}\cdots a_{e_{2k-1}}
&=&(-1)^{\bd(C)-1}a_{e_1}^+\cdots a_{e_{2k-1}}^+
=(-1)^{2k-\fd(C)-1}a_{e_1}^+\cdots a_{e_{2k-1}}^+
\vspace{3pt}\\
&=&(-1)^{\fd(C)+1}a_{e_1}^+\cdots a_{e_{2k-1}}^+
=\sgn(C)w_A(C).
\end{array}
$$
\end{proof}

\begin{proof}[Proof of Theorem \ref{thm:pb}]
The proof proceeds by induction on $n$. 
Let $A\in\mat(2n,2n)$ be skew symmetric with its block decomposition given as above. 
Let $\bar{p}_A=\pairs{p_n,p_{n-1},\ldots,p_0}$ and $\bar{p}_M=\pairs{q_{n-1},q_{n-2},\ldots,q_0}$ 
be Pfaffian sequences for $A$ and $M$ respectively. 
By the inductive hypothesis, we have 
\begin{equation}\label{eq:IH}
q_{n-1-k}=\sum_{\substack{\bar{C}\mbox{ : pclow seq. }\\|\bar{C}|=2k}}\sgn(\bar(C))w_M(\bar{C}). 
\end{equation}
for $0\leq k\leq n-1$. Let $1\leq k\leq n-1$. Then by PB algorithm we have 
\begin{equation}
p_{n-k}=RJ(MJ)^{k-2}S+\sum_{j=3}^kRJ(MJ)^{k-j}Sq_{n-j+1}+a_{12}q_1+q_0
\end{equation}
and 
\begin{equation}
p_0=RJ(MJ)^{n-2}S+\sum_{j=3}^nRJ(MJ)^{n-j}Sq_{n-j+1}+a_{12}q_0
\end{equation}

By equations with the equations (\ref{eq:pclow}) from 
Lemma \ref{lemma:plow} and (\ref{eq:IH}), we get 
\begin{equation}
\begin{aligned}
&RJ(MJ)^{k-j}Sq_{n-j+1}\\
&=\left(\sum_{\substack{C\mbox{ : pclow }\\|\bar{C}|=2(k-j+1),\head(C)=1}}\sgn(C)w_A(C)\right)\cdot
\left(\sum_{\substack{\bar{C}\mbox{ : pclow seq. }\\|\bar{C}|=2(j-1)}}\sgn(\bar{C})w_M(\bar{C})\right)\\
&=\sum_{\substack{C\mbox{ : pclow }\\|\bar{C}|=2(k-j+1),\head(C)=1}}
\sum_{\substack{\bar{C}\mbox{ : pclow seq. }\\|\bar{C}|=2(j-1)}}
\sgn(C)w_A(C)\sgn(\bar{C})w_M(\bar{C})\\
&=\sum_{\substack{C\mbox{ : pclow }\\|\bar{C}|=2(k-j+1),\head(C)=1}}
\sum_{\substack{\bar{C}\mbox{ : pclow seq. }\\|\bar{C}|=2(j-1)}}
\sgn(\pairs{C,\bar{C}})w_A(\pairs{C,\bar{C}})\\
\end{aligned}
\end{equation}

\end{proof}

Note that PB algorithm is a $\sharpl$ algorithm and hence we have 
\begin{corollary}
Pfaffian $\pfaff(A)$ is $\sbi{1}$ definable in $\vsl$. 
\end{corollary}

\section{The proof complexity of Pfaffian}

Some of Pfaffian properties are derivable solely from Pfaffian Berkowitz algorithm. 
Here we present two of them. 

\begin{lemma}($\vsl$)
Let $A\in\mat(2n,2n)$ be skew symmetric and $\lambda$ be any number. Then 
\begin{equation}
\pfaff(\lambda A)=\lambda^n\pfaff(A). 
\end{equation}
\end{lemma}

\begin{theorem}
Let $A\in\mat(2n,2n)$ be skew symmetric, $\vec{q_A}=(q_n,q_{n-1},\ldots,q_0)$ and 
$\vec{r_A}=(r_n,r_{n-1},\ldots,r_0)$ be Pfaffian coefficients of $A$ and $^t\!A$ respectively. 
Then for $0\leq k\leq n$, 
\begin{equation}
r_{n-k}=(-1)^kq_{n-k}.
\end{equation}
Thus $\pfaff(^t\!A)=(-1)^n\pfaff(A)$.
\end{theorem}

Since Pfaffian is a generalization of the determinant, most properties of the determinant 
are given for Pfaffian as well. The difference is that operations on rows or columns on 
$\det(A)$ correspond to operations simultaneously on rows and columns. 

Let $A\in\mat(2n,2n)$ be skew symmetric. Define the following operations:
\begin{itemize}
\item $A[i:j]$ is given by simultaneously swapping rows $i,j$ and swapping columns $i,j$. 
\item $A\pairs{i,j}$ is given by removing rows $i,j$ and columns $i,j$. 
\end{itemize}

Then we have the following properties in analogy with the determinant:

\begin{theorem}[Pfaffian Cofactor Expansion]
Let $A\in\mat(2n,2n)$ be skew symmetric and $1\leq i\leq 2n$. Then 
\begin{equation*}
(\mbox{PCE})\quad 
\pfaff(A)=\sum_{1\leq j\neq i\leq 2n}(-1)^{i+j+\Theta(j-i)}a_{ij}\pfaff(A\pairs{i,j})
\end{equation*}
where $\Theta(k)$ is Heaviside step function. 
\end{theorem}

If we define the determinant by the equation (\ref{eq:det-pf}) then 
properties of the determinant are provable from the corresponding properties 
for Pfaffian in $\vsl$. For instance, we have 

\begin{lemma}[$\vsl$]
(PCE) implies the cofactor expansion of the determinant. 
\end{lemma}

\begin{proof}[Proof Sketch]
The proof is by induction on the number of rows. 
Let $A\in\mat(n,n)$ and 
$B=
\begin{pmatrix}
0&A\\
-^t\!A&0
\end{pmatrix}
$. By applying (PCE) to $B$ and using the inductive hypothesis yields that
\begin{equation*}
\pfaff(B)=
(-1)^{n(n-1)/2}\sum_{1\leq j\leq n}(-1)^{i+j}\det(A_{i,j}).
\end{equation*}
\end{proof}

\begin{theorem}[$\vsl$]
(PCE) implies Cayley's theorem:
\begin{equation}
\all A\in\mat(2n,2n)\mbox{ : skew symmetric}\det(A)=\pfaff(A)^2.
\end{equation}
\end{theorem}

See \cite{okada} for the proof. 

The axiomatic definition of the determinant refers to the multilinearity, 
the alternation and the equation $\det(I)=1$. Similarly, 
the axiomatic definition of Pfaffian (PAD) is the collection of 
the following three statements:
\begin{description}
\item[Multilinearity : ]Let $A(\lambda,i)$ be the matrix $A$ with 
the row and the column $i$ multiplied by $\lambda$. Then $\pfaff(A(\lambda,i))=\lambda\pfaff(A)$. 
\item[Alternation : ]$\pfaff(A[i:j])=-\pfaff(A)$. 
\item[Identity : ]$\pfaff(J)=-1$.
\end{description}

\begin{theorem}\label{thm:mult-id}
$\vsl$ proves Multilinearity on the first row and column and Identity.
\end{theorem}

\begin{proof}
The first part is easy. For the second part, let $J_n\in\mat(2n,2n)$. 
Then 
\begin{equation*}
P_{J_n}=
\begin{pmatrix}
1&0&\ddots&\\
1&1&\ddots&\\
\vdots&\vdots&\ddots&1\\
0&0&\cdots&1
\end{pmatrix}
\end{equation*}
From this we have the recurrence $\pfaff(J_n)=\pfaff(J_{n-1})$. 
Since we have $\pfaff(J_1)=-1$, the claim is immediate. 
\end{proof}

For the determinant, cofactor expansion and the axiomatic definition are equivalent 
in $\vsl$. This is also the case for Pfaffian.

\begin{theorem}[$\vsl$]
(PCE) and (PAD) are equivalent. 
\end{theorem}

\begin{proof}
First we show that (PALT) implies (PCE). 
Let $A\in\mat(2n,2n)$ be skew symmetric and $1<i\leq 2n$. 
Apply (PALT) for rows and columns $1,i$ yields $\pfaff(A\pairs{1,i})=-\pfaff(A)$. 
By Theorem \ref{thm:mult-id}, we can expand $A\pairs{i,j}$ on the first row and column. 
Then applying (PALT) again yields (PCE). 

For the other direction, we can show that (PCE) implies (PALT). 

\end{proof}

\begin{theorem}[$\vsl$]\label{thm:det(I)=1}
Let $I\in\mat(2n,2n)$ be the identity matrix. Then
\begin{equation}
\pfaff
\begin{pmatrix}
0&I\\
-I&0
\end{pmatrix}
=(-1)^n.
\end{equation}
\end{theorem}

\begin{proof}
Let 
$C_0=
\pfaff
\begin{pmatrix}
0&I\\
-I&0
\end{pmatrix}
\in\mat(4n,4n)$ and 
\begin{equation*}
C_k=
\begin{pmatrix}
0&0&R_{k+1}\\
0&0&-^t\!S_{k+1}\\
-^t\!R_{k+1}&S_{k+1}&C_{k+1}\\
\end{pmatrix}
\in\mat(4n-2k,4n-2k). 
\end{equation*}
Then $C_k$ is of the form
\begin{equation*}
\begin{pmatrix}
&&I_{2n-2k}\\
&O_{2k}&\\
-I_{2n-2k}&&\\
\end{pmatrix}
\in\mat(4n-2k,4n-2k)
\end{equation*}
where all blank entries are zero. 

Let 
$\vec{q_k}=
\begin{pmatrix}
q^k_{2n-k},q^k_{2n-k-1},\ldots,q^k_0
\end{pmatrix}
$ be Berkowitz sequence for $C_k$. From Berkowitz algorithm, it follows that 
\begin{equation*}
\begin{pmatrix}
q^k_{2n-k}
\vspace{3pt}\\
q^k_{2n-k-1}
\vspace{3pt}\\
\vdots
\vspace{3pt}\\
q^k_0
\end{pmatrix}
=
\begin{pmatrix}
1&0&&&&\\
0&1&\ddots&&&\\
-1&0&\ddots&&&\\
\vdots&\vdots&\ddots&&\\
&&&&1&\\
&&&&0&1\\
0&0&0&\cdots&-1&0
\end{pmatrix}
\begin{pmatrix}
q^{k+1}_{2n-k-1}
\vspace{3pt}\\
q^{k+1}_{2n-k-2}
\vspace{3pt}\\
\vdots
\vspace{3pt}\\
q^{k+1}_0
\end{pmatrix}
\end{equation*}

Hence we have the following recurrences:
\begin{equation*}
\begin{aligned}
&q^k_{n-k}=1,\ q^k_{n-k-1}=q^{k+1}_{n-k-2},\\
&q^k_{n-k-i}=-q^{k+1}_{n-k-i+1}+q^{k+1}_{n-k-i-1}\ (2\leq i<2n-k),\\
&q^k_0=-q^k_1.
\end{aligned}
\end{equation*}

We claim that $\all j<k\;p^k_j=0$. This is proved by backward induction on $k\leq n$. 
If $k=n$ then 
$\vec{q_n}=
\begin{pmatrix}
1&0&\cdots&0
\end{pmatrix}
$ so the claim is obvious. Suppose that $\all j<k\;p^k_j=0$. and $j\leq k$. Then 
by the inductive hypothesis, 
\begin{equation*}
\begin{aligned}
&q^k_{2n-k-1}=q^{k+1}_{2n-k-2}=0,\\
&q^k_{n-k-i}=-q^{k+1}_{n-k-i+1}+q^{k+1}_{n-k-i-1}=0. 
\end{aligned}
\end{equation*}

Hence we have $q^k_k=-q^{k+1}_{k+1}$ and by backward induction, we conclude that 
$\pfaff(C_0)=(-1)^n$. 
\end{proof}

The following equation can be regarded as the analogue of the multiplicativity. 

\begin{theorem}
Let $A\in\mat(2n,2n)$ be skew symmetric and $B\in\mat(2n,2n)$. Then 
\begin{equation*}
(MP)\quad \pfaff(^t\!BAB)=\pfaff(A)\det(B).
\end{equation*}
\end{theorem}

Note that if $B$ is skew symmetric then by Cayley's theorem (equation \ref{eq:cayley}), 
we have $\det(B)=\pfaff(B)^2$. 
Hence we obtain
$$
\pfaff(^t\!BAB)=\pfaff(A)\pfaff(B)^2.
$$

We expect that (MP) implies most properties of Pfaffian in $\vsl$. 
The rest of this section is devoted to the consideration of this problem.

First remark that the condition of the axiomatic definition other than Alternation 
are provable in $\vsl$. 

\begin{theorem}\label{thm:multi}
$\vsl$ proves Multilinearity for $i=1$ and Identity.
\end{theorem}

\begin{proof}
Multilinearity for the first row and column is straightforward from PB algorithm. 
Identity is proved by induction on $n$. 
\end{proof}

Suppose that Alternation is proved form (MP). Then it is easy to see that (MP) implies 
other properties of Pfaffian, namely (PCE) and Cayley's Theorem. 
To prove Alternation from (MP), remark that 
$$
A[i:j]=^t\!I_{ij}AI_{ij}. 
$$
So by (MP) we have 
$$
\pfaff(A[i:j])=\pfaff(^t\!I_{ij}AI_{ij})=\pfaff(A)\det(I_{ij}). 
$$
Hence it suffices to show that 
\begin{equation}\label{eq:identity}
\det(I_{ij})=-\det(I)
\end{equation} 
is provable for the identity matrix $I$ of even order. However, it seems fairly complicated to 
directly prove the equation (\ref{eq:identity}). So we argue in a simpler manner. 

First note that
\begin{theorem}[$\vsl$]\label{thm:det(I[k])=-1}
Let $I\in\mat(2n,2n)$ be the identity matrix and for $1\leq k<2n$, 
$I[k]$ be the alternation of rows $k$ and $k+1$ in $I$. Then 
\begin{equation}
\pfaff
\begin{pmatrix}
0&I[k]\\
-I[k]&0
\end{pmatrix}
=(-1)^{n+1}. 
\end{equation}
\end{theorem}

\begin{proof}
We argue similarly as in Theorem \ref{thm:det(I)=1}. The proof is divided into two cases. 

\vspace{6pt}
\noindent
Case 1: $k\mod 2=1$. Let $k=2l+1$ and 
\begin{equation*}
I[k]=
\begin{pmatrix}
I_{2l}&&\\
&&1&\\
&1&&\\
&&&I_{2m}
\end{pmatrix}
\end{equation*}
where $n=l+m+1$. Define 
\begin{equation*}
C_0=
\begin{pmatrix}
0&I[k]\\
-I[k]&0
\end{pmatrix},\ 
C_{j}=
\begin{pmatrix}
0&0&R_{j+1}\\
0&0&^t\!S_{j+1}\\
^t\!R_{j+1}&S_{j+1}&C_{j+1}
\end{pmatrix}\ 
(0\leq k<n). 
\end{equation*}
Then for $0\leq j<l$, $j=l$ and $l<j<n$, $C_j$ is of the form 
\begin{equation*}
\begin{aligned}
&\begin{pmatrix}
&&I[2(l-j)+1]\\
&O_{2j}&\\
-I[2(l-j)+1]&&\\
\end{pmatrix}
\in\mat(4n-2j,4n-2j),\\
&\begin{pmatrix}
&&&&&1&&\\
&&&&1&&\\
&&&&&&I\\
&&&O_{2l}&&&\\
&-1&&&&\\
-1&&&&&\\
&&-I&&&\\
\end{pmatrix}
\in\mat(4n-2l,4n-2l), \mbox{and}\\
&\begin{pmatrix}
&&I\\
&O_{2j}&\\
-I&&\\
\end{pmatrix}
\in\mat(4n-2j,4n-2j),\\
\end{aligned}
\end{equation*}
respectively. Moreover, $C_n=0$. 

Remark that Berkowitz matrix for $C_j$ with $j\neq l$ are the same as in Theorem \ref{thm:det(I)=1}. 
For $C_l$, we have 
\begin{equation*}
P_{C_l}=
\begin{pmatrix}
1&0&&&&\\
0&1&\ddots&&&\\
1&0&\ddots&&&\\
\vdots&\vdots&\ddots&&\\
&&&&1&\\
&&&&0&1\\
0&0&0&\cdots&1&0
\end{pmatrix}
\end{equation*}

Now the proof is identical to that for Theorem \ref{thm:det(I)=1}. Specifically, we have 
\begin{equation*}
q_j^j=
\begin{cases}
-q_{j+1}^{j+1}&\mbox{ if }0\leq j\leq l<n, \\
q_{j+1}^{j+1}&\mbox{ if }j=l. 
\end{cases}
\end{equation*}
and from this recurrence, the claim follows immediately. 

\vspace{6pt}
\noindent
Case 1: $k\mod 2=0$. Let $k=2l$. By a similar block decomposition as in Case 1, we have 
\begin{equation*}
C_{l-1}=
\begin{pmatrix}
&&B\\
&O_{2(l-1)}&\\
-B&&
\end{pmatrix},\mbox{ and }
C_l=
\begin{pmatrix}
&&B'\\
&O_{2l-1}&\\
-B'&&
\end{pmatrix}
\end{equation*}
where 
\begin{equation*}
B=
\begin{pmatrix}
1&&&&\\
&&1&&\\
&1&&&\\
&&&1&\\
&&&&I\\
\end{pmatrix},\mbox{ and }
B'=
\begin{pmatrix}
1&&&\\
&&1&\\
&&&I\\
\end{pmatrix}.
\end{equation*}

Berkowitz matrices for these two matrices are 
\begin{equation*}
P_{C_{l-1}}=
\begin{pmatrix}
1&0&&&&\\
0&1&\ddots&&&\\
0&0&\ddots&&&\\
-1&0&\ddots&&&\\
\vdots&\vdots&\ddots&&\\
&&&&1&\\
&&&&0&1\\
0&0&0&\cdots&0&0
\end{pmatrix}\mbox{ and }
P_{C_l}=
\begin{pmatrix}
1&0&&&&\\
0&1&\ddots&&&\\
0&0&\ddots&&&\\
\vdots&\vdots&\ddots&&\\
&&&&1&\\
&&&&0&1\\
0&0&0&\cdots&0&0
\end{pmatrix}. 
\end{equation*}

From these computations we obtain a recurrence which implies the claim. 
\end{proof}

\begin{theorem}[$\vsl$]\label{thm:(MP)->(PALT)}
(MP) implies (PALT). 
\end{theorem}

\begin{proof}
Let $A\in\mat(2n,2n)$ be skew symmetric and $1\leq i\neq j\leq 2n$. Then we can 
effectively construct a sequence $1\leq k_1,k_2,\ldots,k_{2l-1}\leq 2n$ such that 
\begin{equation*}
A[i:j]={^t\!I_{k_{2l-1}}}\cdots{^t\!I_{k_2}}{^t\!I_{k_1}}AI_{k_1}I_{k_2}\cdots I_{k_{2l-1}}. 
\end{equation*}

By applying (MP) repeatedly to $A[i:j]$, we obtain 
\begin{equation*}
\pfaff(A[i:j])=\pfaff(A)\det(I_{k_1})\det(I_{k_2})\cdots\det(I_{k_{2l-1}})
\end{equation*}
Note that this is where we require induction in $\vsl$. Now the claim follows immediately 
from Theorem \ref{thm:det(I[k])=-1}. 
\end{proof}

\section{Cayley-Hamilton Theorem for Pfaffian}

So far we have seen that Pfaffian and the determinant have a lot of common properties. 
Hence one might ask whether Cayley-Hamilton type theorem is possible for Pfaffian. 
The answer is yes and in this section we present a version of Cayley-Hamilton theorem for Pfaffian. 

\begin{definition}
Let $A\in\mat(2n,2n)$ be skew symmetric and $q^A_n,q^a_{n-1},\ldots,q^A_0$ be 
its Pfaffian coefficients. Define Pfaffian characteristic polynomial as
\begin{equation}
\Phi_A(x)=q^A_nx^n+q^A_{n-1}x^{n-1}+\cdots+q^A_0.
\end{equation}
\end{definition}

\begin{theorem}[Pfaffian Cayley-Hamilton Theorem]\label{thm:pch}
Let $A\in\mat(2n,2n)$ be skew symmetric. Then 
\begin{equation}
\Phi_A(AJ)=q^A_n(AJ)^n+q^A_{n-1}(AJ)^{n-1}+\cdots+q^A_0I=0.
\end{equation}
\end{theorem}

We can prove this theorem in several ways. 
One way is to use the combinatorial argument which is used to prove Cayley-Hamilton Theorem 
for the determinant due to Straubing \cite{straubing}. However, we do not know whether such 
proof can be formalized in $\vsl$. 

Here we give a proof from (PCE) which can be formalized in $\vsl$. 

\begin{proof}[Proof of Theorem \ref{thm:pch}]
Let $A\in\mat(2n,2n)$ be skew symmetric and for $1\leq i\neq j\leq 2n$. Define 
\begin{equation}\label{eq:aij}
\tilde{a}_{ij}=(-1)^{i+j+\Theta(j-i)}\pfaff(A\pairs{i,j})
\end{equation}
and
\begin{equation*}
A^*=
\begin{pmatrix}
0&\tilde{a}_{21}&\cdots&\tilde{a}_{n1}\\
\tilde{a}_{12}&0&\cdots&\tilde{a}_{n2}\\
\vdots&\vdots&\ddots&\vdots\\
\tilde{a}_{1n}&0&\cdots&\tilde{a}_{nn}\\
\end{pmatrix}\in\mat(2n,2n).
\end{equation*}

By (PCE), we have 
\begin{equation}\label{eq:A*}
AA^*=\pfaff(A)I.
\end{equation}

On the other hand, let 
\begin{equation*}
\Padj(A)=-J((AJ)^{n-1}+p^A_{n-1}(AJ)^{n-2}+\cdots+p^A_1I). 
\end{equation*}
Then (PCH) can be expressed as 
\begin{equation}\label{eq:padj}
A\Padj(A)=\pfaff(A)I
\end{equation}

For $1\leq i\neq j\leq 2n$, let
\begin{equation*}
C_{ij}=
\begin{pmatrix}
0&0&e_i\\
0&0&-^t\!f_j\\
-^t\!e_i&f_j&A
\end{pmatrix}
\end{equation*}
where $e_i$ is the row vector with all entries $0$ except 
the $i$th entry $1$ and $f_j$ is the column vector with all entries $0$ except 
the $j$th entry $-1$. By definition, we have 
\begin{equation}\label{eq:cij}
\pfaff(C_{ij})=c_{12}\pfaff(A)-e_i\Padj(A)f_j=\Padj(A)_{ij}
\end{equation}

From (\ref{eq:aij}), (\ref{eq:A*}), (\ref{eq:padj}) and (\ref{eq:cij}), it suffices to show that 
\begin{equation}
\pfaff(C_{ij})=(-1)^{i+j+\Theta(i-j)}\pfaff(A\pairs{i,j}). 
\end{equation}
The proof is divided into two case. 

First suppose that $i<j$. We will show that 
$\pfaff(C_{ij})=(-1)^{i+j}\pfaff(A\pairs{i,j})$. By alternation, we transform $C_{ij}$ 
into the form 
\begin{equation}
C'_{ij}=
\begin{pmatrix}
0&1&0&\cdots&0\\
-1&0&*&\cdots&*\\
0&0&&&\\
\vdots&\vdots&&B&\\
0&0&&&
\end{pmatrix}, 
B=
\begin{pmatrix}
0&f_{j-1}\\
f_{j-1}&A\pairs{i,-}
\end{pmatrix}
\end{equation}
where $A\pairs{i,-}$ is obtained from $A$ by removing $i$th row and column. 
This transformation requires $i$ alternations. So we have $\pfaff(C'_{ij})=(-1)^i\pfaff(C_{ij})$. 
By Berkowitz algorithm $\pfaff(C'_{ij})=\pfaff(B)$ and applying (PCE) to the first row of $B$ yields 
$\pfaff(B)=(-1)^{1+j+1}\pfaff(A\pairs{i,j})$. Thus we obtain
\begin{equation}
\pfaff(C_{ij})=(-1)^i(-1)^{1+j+1}\pfaff(A\pairs{i,j})=(-1)^{i+j}\pfaff(A\pairs{i,j}). 
\end{equation}

Second, suppose that $j<i$. We want show that 
$\pfaff(C_{ij})=(-1)^{i+j+1}\pfaff(A\pairs{i,j})$. 
The proof is almost identical to the case for $i<j$. By alternation 
we transform $C_{ij}$ into $C'_{ij}$ where 
\begin{equation}
B=
\begin{pmatrix}
0&f_j\\
f_j&A\pairs{i,-}
\end{pmatrix}
\end{equation}
This is because $j$th row and column shift by $1$ by applying $i$ alternations to $C_{ij}$. 
Then the rest of the proof is exactly the same as before. 
\end{proof}

The converse is also true, that is, 
\begin{theorem}[$\vsl$]
(PCH) implies (PCE). Thus (PCE), (PCH) and (PAD) are equivalent over $\vsl$. 
\end{theorem}

We omit the proof due to the limitation of space. 

\section{Berkowitz algorithm for Pfaffian pairs}

In \cite{rote}, it is shown that the product of pfaffians are computable by means of alternating 
clow sequences. This fact leads to a Berkowitz type algorithm computing $\pfaff(A)\pfaff(B)$. 
In this section we present such an algorithm. 

Let $1\leq i\leq n$ and $k$ be a number. We define

\begin{equation}
\begin{aligned}
&\mathcal{C}_{i,k}
=\left\{\bar{C} : \mbox{alternating clow},\ \head(\bar{C})\geq i,\ |\bar{C}|=2k\right\},\mbox{ and }\\
&\mathcal{D}_{i,k}
=\left\{\bar{C} : \mbox{alternating clow},\ \head(\bar{C})=i,\ |\bar{C}|=2k\right\}.\\
\end{aligned}
\end{equation}

For skew symmetric matrices $A,B\in\mat(n,n)$ let $m=\lfloor n/2\rfloor$. 
We define P-coefficients of $A$ and $B$ as 
\begin{equation}\label{eq:p-coeff}
\vec{q}_{A,B}=(q_m,q_{m-1},\ldots,q_0)
\end{equation}
where
\begin{equation}\label{eq:p-coeff2}
\begin{aligned}
&q_m=1,\\
&q_{m-k}=\sum_{\bar{C}\in\mathcal{C}_{1,k}}\sgn(\bar{C})w_{A,B}(\bar{C}),\mbox{ for }1\leq k<n,\mbox{ and}\\
&q_0=
\left\{
\begin{aligned}
\sum_{\bar{C}\in\mathcal{C}_{1,n}}\sgn(\bar{C})w_{A,B}(\bar{C})&\mbox{ if }n\equiv 0\pmod 1,\\
\sum_{\bar{C}\in\mathcal{D}_{1,n}}\sgn(\bar{C})w_{A,B}(\bar{C})&\mbox{ if }n\equiv 0\pmod 0.
\end{aligned}
\right.
\end{aligned}
\end{equation}
where for a clow $C=\pairs{e_1,e_2,\ldots,e_{2k-1},e_{2k}}$, we define the weight 
\begin{equation}
w_{A,B}(C)=a_{e_1}b_{e_2}\cdots a_{e_{2k-1}}b_{e_{2k}}
\end{equation}
and for a clow sequence $\bar{C}=\pairs{C_1,\ldots,C_l}$
\begin{equation}
w_{A,B}(\bar{C})=\prod_{1\leq i\leq l}w_{A,B}(C_i). 
\end{equation}

In \cite{rote}, it is shown that $\pfaff(A)\pfaff(B)=q_0$. So this notion is a generalization 
of the clow presentation of Pfaffian pairs. 

We will construct a recursive algorithm which computes P-coefficients (\ref{eq:p-coeff}). 
Let 
\begin{equation*}
A=
\begin{pmatrix}
0&R\\
-^t\!R&M
\end{pmatrix},\ 
B=
\begin{pmatrix}
0&-^tS\\
S&N
\end{pmatrix}
\end{equation*}

Let $1\leq k<n$. By the equation (\ref{eq:p-coeff}), we have
\begin{equation}
q_{m-k}
=\sum_{\bar{C}\in\mathcal{D}_{1,k}}\sgn(\bar{C})w_{A,B}(\bar{C})+
\sum_{\bar{C}\in\mathcal{C}_{1,k}}\sgn(\bar{C})w_{A,B}(\bar{C}).
\end{equation}

Let 
\begin{equation}
c(l,k,A,B)=\sum_{\bar{C}\in\mathcal{C}_{1,k}}w_{A,B}(\bar{C}),\\
d(l,k,A,B)=\sum_{\bar{C}\in\mathcal{D}_{1,k}}w_{A,B}(\bar{C}).
\end{equation}
Then for $2\leq l\leq k$, 
\begin{equation}
\sum_{\bar{C}\in\mathcal{D}_{1,k}}\sgn(\bar{C})w_{A,B}(\bar{C})
=\sum_{1\leq j\leq k-l+1}d(1,j,A,B)c(l-1,k-1,M,N). 
\end{equation}
So we have 
\begin{equation}\label{eq:q_{m-k}}
\begin{aligned}
&\sum_{\bar{C}\in\mathcal{D}_{1,k}}\sgn(\bar{C})w_{A,B}(\bar{C})\\
&=-d(1,k,A,B)+\sum_{2\leq l\leq k}\sum_{1\leq j\leq k-l+1}(-d(1,j,A,B))(-1)^{l-1}c(l-1,k-j,M,N)\\
&=-d(1,k,A,B)+\sum_{1\leq j\leq k-1}(-d(1,j,A,B))\sum_{1\leq l\leq k-j}(-1)^{l}c(l,k-j,M,N)\\
\end{aligned}
\end{equation}
Remark that 
\begin{equation*}
d(1,j,A,B)=R(NM)^{j-1}S
\end{equation*}
and
\begin{equation*}
\sum_{1\leq l\leq k}(-1)^{l}c(l,k,M,N)
=
\left\{
\begin{aligned}
&r_{m-k+j}&\mbox{ if }n\equiv 1\pmod 2\\
&r_{m-k+j-1}&\mbox{ if }n\equiv 0\pmod 2\\
\end{aligned}
\right.
\end{equation*}
By substituting these equations in (\ref{eq:q_{m-k}}) we obtain the following:

\begin{theorem}
Let $A,B\in\mat(n,n)$, $m$, $\vec{q}_{A,B}$ and $\vec{q}_{M,N}$ be as above. Then 
for $1\leq k<n$, 
\begin{equation}
q_{m-k}=
\left\{
\begin{aligned}
&-R(NM)^{k-1}S-\sum_{1\leq j\leq k-1}R(NM)^{j-1}S\;r_{m-k+j}+r_{m-k}&
\mbox{ if }n\equiv 1\pmod 2,\\
&-R(NM)^{k-1}S-\sum_{1\leq j\leq k-1}R(NM)^{j-1}S\;r_{m-k+j-1}+r_{m-k-1}&
\mbox{ if }n\equiv 0\pmod 2,\\
\end{aligned}
\right.
\end{equation}
and 
\begin{equation}
q_0=
\left\{
\begin{aligned}
&-R(NM)^{n-1}S-\sum_{1\leq j\leq k-1}R(NM)^{j-1}S\;r_{m-k+j}+r_0&
\mbox{ if }n\equiv 1\pmod 2,\\
&-R(NM)^{n-1}S-\sum_{1\leq j\leq k-1}R(NM)^{j-1}S\;r_{m-k+j-1}&
\mbox{ if }n\equiv 0\pmod 2,\\
\end{aligned}
\right.
\end{equation}
\end{theorem}
Now we are ready to construct a recursive procedure which computes $\vec{q}_{A,B}$ from 
$\vec{q}_{M,N}$. The procedure is divided into two cases. 

\begin{theorem}\label{thm:pb pairs}
Let $A,B\in\mat(n,n)$, $m$, $\vec{q}_{A,B}$ be as above. 
Let 
\begin{equation*}
\vec{q}_{M,N}=
\left\{
\begin{aligned}
&(r_{m-1},r_{m-2},\ldots,r_0)&\mbox{ if }n\equiv 0\pmod 2\\
&(r_{m},r_{m-1},\ldots,r_0)&\mbox{ if }n\equiv 1\pmod 2.\\
\end{aligned}
\right.
\end{equation*}
If $n\equiv 0\pmod 2$ then 
\begin{equation}
\begin{pmatrix}
q_m\\
q_{m-1}\\
q_{m-2}\\
\vdots\\
q_0
\end{pmatrix}
=
\begin{pmatrix}
1&0&\\
-RS&1&\ddots\\
-R(NM)S&-RS&\ddots\\
\vdots&\vdots&\ddots\\
-R(NM)^{m-1}S&-R(NM)^{m-2}S&\cdots&-RS
\end{pmatrix}
\begin{pmatrix}
r_{m-1}\\
r_{m-2}\\
\vdots\\
r_0
\end{pmatrix}
\end{equation}
and if $n\equiv 1\pmod 2$ then 
\begin{equation}
\begin{pmatrix}
q_m\\
q_{m-1}\\
q_{m-2}\\
\vdots\\
q_0
\end{pmatrix}
=
\begin{pmatrix}
1&0&0\\
-RS&1&0&\ddots\\
-R(NM)S&-RS&\ddots\\
\vdots&\vdots&\ddots&\ddots\\
-R(NM)^{m-1}S&-R(NM)^{m-2}S&\cdots&-RS&1
\end{pmatrix}
\begin{pmatrix}
r_m\\
r_{m-1}\\
r_{m-2}\\
\vdots\\
r_0
\end{pmatrix}
\end{equation}
\end{theorem}

Hence we have two ways to compute Pfaffian in $\vsl$; one by PB algorithm in 
Definition \ref{def:pb algorithm} and the other by computing $\pfaff(A)\pfaff(J)$ 
by the algorithm for Pfaffian pairs. Moreover, $\vsl$ proves that these two 
definitions of Pfaffians coincide. 

\begin{theorem}[$\vsl$]
Let $A\in\mat(2n,2n)$ be skew symmetric, $\vec{q}_{A}$ be P-coefficients of $A$ 
and $\vec{q}_{A,J}$ be coefficients of P-coefficients of $A$, $J$. Then 
$\vec{q}_{A}=\vec{q}_{A,J}$.
\end{theorem}

\begin{proof}
By induction on $n$. Let $A,J\in\mat(2n,2n)$ and define block decompositions as 
\begin{equation*}
A=
\begin{pmatrix}
0&R^+\\
-^t\!R^+&M
\end{pmatrix},\ 
M=
\begin{pmatrix}
0&-^t\!S\\
S&N
\end{pmatrix},\ 
J=
\begin{pmatrix}
0&-^t\!P\\
P&K
\end{pmatrix},\ 
K=
\begin{pmatrix}
0&-^t\!Q\\
Q&L
\end{pmatrix}
\end{equation*}
where 
$R^+=
\begin{pmatrix}
a_{12}&R
\end{pmatrix}$. Note that 
\begin{equation*}
P=
\begin{pmatrix}
-1\\
0\\
\vdots\\
0
\end{pmatrix},\ 
Q=
\begin{pmatrix}
0\\
\vdots\\
0
\end{pmatrix}.
\end{equation*}

So we have 
\begin{equation*}
R^+P=
\begin{pmatrix}
a_{12}&R
\end{pmatrix}
\begin{pmatrix}
-1\\
0
\end{pmatrix}=-a_{12}
\end{equation*}
and 
\begin{equation*}
KM=
\begin{pmatrix}
0&-^t\!Q\\
Q&L
\end{pmatrix}
\begin{pmatrix}
0&-^t\!S\\
S&N
\end{pmatrix}
=
\begin{pmatrix}
-^t\!QS&-^t\!QN\\
LS&-^t\!QS+LN
\end{pmatrix}
=
\begin{pmatrix}
0&0\\
LS&LN
\end{pmatrix}. 
\end{equation*}
Therefore we obtain
\begin{equation*}
R^+(KM)^kP=
\begin{pmatrix}
a_{12}&R
\end{pmatrix}
\begin{pmatrix}
0&0\\
(LN)^{k-1}LS&(LN)^k
\end{pmatrix}
\begin{pmatrix}
-1\\
0
\end{pmatrix}
=-R(LN)^{k-1}LS
\end{equation*}

Let $\vec{q}_{A,J}$, $\vec{q}_{M,K}$ and $\vec{q}_{N,L}$ be P-coefficients of 
$A,J$, $M,K$ and $N,L$ respectively. By Theorem \ref{thm:pb pairs} we have 
\begin{equation}\label{eq:pb A,J}
\vec{q}_{A,J}=
\begin{pmatrix}
1&0&\\
a_{12}&1&\ddots\\
RLS&a_{12}&\ddots\\
\vdots&\vdots&\ddots\\
R(LN)^{k-2}LS&R(LN)^{k-3}LS&\cdots&a_{12}\\
\end{pmatrix}
\vec{q}_{M,K}.
\end{equation}

Since $Q=0$ we have 
\begin{equation}
\vec{q}_{M,K}=I\vec{q}_{N,L}. 
\end{equation}

Now by the inductive hypothesis, $\vec{q}_{N,L}=\vec{q}_N$. 
From this and the equation (\ref{eq:pb A,J}) we obtain the claim. 
\end{proof}

\section{Closing remarks}

We have proved that basic Pfaffian properties are provable from Pfaffian version of 
multiplicativity (MP). This fact is similar to that for the determinant where 
(MP) is replaced by the multiplicativity of the determinant. 
So it is natural to conjecture that (MP) is provable in $\vnc{2}$. 

In $\vsl$, the determinant can be defined in two ways; one by Berkowitz algorithm 
for the determinant and the other defined from Pfaffian by the equation (\ref{eq:det-pf}).
It is easily seen that these two definitions are equivalent if we admit (PCE). 

A more challenging problem is to prove the multiplicativity of the determinant
\begin{equation}\label{eq:md}
\det(AB)=\det(A)\det(B)
\end{equation}
in $\vsl$ or some extension of it which is a subsystem of $\vnc{2}$. 
The algorithm given in last section might be a first step towards this 
since both sides of (\ref{eq:md}) can be expressed by way of pclow sequences. 

The ultimate goal of this work is to give candidate hard tautologies for Frege 
proof system which have quasi-polynomial Frege proofs. 
Since such tautologies are provable in weak systems which can formalize Pfaffian 
we believe that such candidate are consequences of Pfaffian identities 
considered in this article. So finding proofs of combinatorial theorems 
from Pfaffian will be our next step.

\end{document}